\documentclass[11pt]{article}
\usepackage{curves}
\usepackage{amsmath, amssymb, lamsarrow}
\usepackage{amssymb}
\def\mineappendix{
        \setcounter{section}{1}
        \setcounter{subsection}{0}
        \def\thesection{\Alph{section}}
        \def\sectionap{\@startsection  {section}{1}{\z@}
                        {-3.5ex plus-1ex minus-.2ex} {0ex plus.2ex}
                        {\reset@font\Large\bf  Appendix:  \, }
                        }
        }
\makeatother
\def\Proclaim #1. #2\par{\bigbreak\noindent{\sc#1.\enspace}{\it#2}\par}

\newcommand{\gwii}[1]{\left< \hspace{-2pt} \left< \, #1 \,
        \right>  \hspace{-2pt} \right>_{0}}

\newcommand{\gwiione}[1]{\left< \hspace{-2pt} \left< \, #1 \,
        \right> \hspace{-2pt} \right>_{1}}

\newcommand{\gwiig}[1]{\left< \hspace{-2pt} \left< \, #1 \,
    \right> \hspace{-2pt} \right>_{g}}
\newcommand{\gwiih}[2]{\left< \hspace{-2pt} \left< \, #2 \,
    \right> \hspace{-2pt} \right>_{#1}}

\newcommand{\ga}{\gamma_{\alpha}}
\newcommand{\gua}{\gamma^{\alpha}}
\newcommand{\gb}{\gamma_{\beta}}
\newcommand{\gub}{\gamma^{\beta}}
\newcommand{\gm}{\gamma_{\mu}}
\newcommand{\gum}{\gamma^{\mu}}

\newcommand{\qp}{\circ}

\newtheorem{lem}{Lemma}[section]

\newtheorem{thm}[lem]{Theorem}

\title{Genus-1 Virasoro conjecture along quantum volume direction}

\author{Xiaobo Liu
\thanks{Research was partially supported by NSF grant
DMS-0905227.} }

\date{}

\begin{document}

\maketitle

\begin{abstract}
In this paper, we show that the derivative of the genus-1 Virasoro
conjecture for Gromov-Witten invariants along the direction of
quantum volume element holds for all smooth projective varieties.
This result provides new evidence for the  Virasoro
conjecture.
\end{abstract}

\section{Introduction}

The Virasoro conjecture predicts that the generating functions of
the Gromov-Witten invariants of smooth projective varieties are
annihilated by a sequence of differential operators which form a
half branch of the Virasoro algebra. This conjecture was proposed by
Eguchi-Hori-Xiong \cite{EHX} and modified by S. Katz \cite{CK}. In
case the underlying manifold is a point, this conjecture is
equivalent to Witten's conjecture \cite{W}, proved by Kontsevich
\cite{K}, that the generating function of intersection numbers on
the moduli spaces of stable curves is a $\tau$-function of the KdV
hierarchy. Together with Tian, we proved that the genus-0 part of
the Virasoro conjecture holds for all compact symplectic manifolds
(cf. \cite{LT}). For manifolds with semisimple quantum cohomology,
the genus-1 part of this conjecture was proved by Dubrovin and Zhang
\cite{DZ}. Without assuming semisimplicity, the genus-1 Virasoro
conjecture was studied in \cite{L1} and \cite{L2}. Among other
results, it was proved in \cite{L1} that the genus-1 Virasoro
conjecture can be reduced to the the $L_1$-constraint. Using the genus-1
topological recursion relation, it was also proved that Virasoro
constraints can be reduced to equations on the {\it small phase
space}, i.e. the space of cohomology classes of the underlying
manifold. Compatibility conditions for Virasoro conjectures were
studied in \cite{L2}. Despite these efforts, the general case of the
genus-1 Virasoro conjecture is still largely open. In this paper, we
give more evidence to the genus-1 Virasoro conjecture without any
assumption on the quantum cohomology of the underlying manifold.

Let $M$ be a smooth projective variety.
Choose a basis $\{ \gamma_{\alpha} \mid \alpha = 1, \ldots, N\}$ of
the space of cohomology  classes $H^{*}(M; \mathbb{C})$. For simplicity, we assume
$H^{\rm odd}(M; \mathbb{C}) = 0$. We choose the basis in such a way that $\gamma_{1}$
is the identity of the cohomology ring and $\gamma_{\alpha} \in
H^{p_{\alpha}, q_{\alpha}}(M)$ for some integers $p_{\alpha}$ and
$q_{\alpha}$. Let $\{ t^{1}, \ldots, t^{N}\}$ be the coordinates on
$H^{*}(M; \mathbb{C})$ with respect to this basis. We can identify
each $\gamma_{\alpha}$ with the vector field
$\frac{\partial}{\partial t^{\alpha}}$ and further identify each
cohomology class with a constant vector field on $H^{*}(M; \mathbb{C})$. Let
\begin{equation} \label{eqn:balpha}
         b_{\alpha} = p_{\alpha} - \frac{1}{2}(d-1)
\end{equation}
where $d$ is the complex dimension of $M$.
 Then the {\it Euler vector field} (on the small phase space)
 is defined to be
 \[ E := c_{1}(M) + \sum_{\alpha} (b_{1} + 1 - b_{\alpha})
            t^{\alpha} \gamma_{\alpha} . \]

We refer to \cite{LiT} \cite{RT} for definitions of Gromov-Witten invariants.
In genus-$1$ case, it suffices to study only primary Gromov-Witten invariants
since all genus-$1$ descendant invariants can be reduced to primary invariants
due to the genus-$1$ topological recursion relation. Therefore we only
consider primary Gromov-Witten invariants in this paper.
Let $F_g$ be the generating function of genus-$g$ primary Gromov-Witten invariants
of $M$. The {\it $k$-point function} is defined to be
\[ \gwiig{v_1 \cdots v_k} := \sum_{\alpha_{1}, \ldots, \alpha_{k}}
        f^{1}_{\alpha_{1}} \cdots f^{k}_{\alpha_{k}} \,\,
        \frac{\partial^k F_g}{\partial t^{\alpha_1} \cdots \partial t^{\alpha_k}},
\]
for vector fields
$v_{i} = \sum_{\alpha}
    f^{i}_{\alpha} \gamma_{\alpha}$ where
$f^{i}_{\alpha}$ are functions on  $H^*(M; \mathbb{C})$.
Note that $F_g$ and $\gwiig{\cdots}$ in
this paper corresponds to $F_g^s$ and $\gwiih{g, s}{\cdots}$
in \cite{L1}.
Let $\eta_{\alpha \beta} = \int_{M} \gamma_{\alpha} \cup \gamma_{\beta}$
be the intersection form on $H^{*}(M, {\Bbb C})$.
We will use $\eta = (\eta_{\alpha \beta})$ and
$\eta^{-1} = (\eta^{\alpha \beta})$ to lower and raise indices.
For example
$\gua := \eta^{\alpha \beta} \gb$ where repeated indices should be summed
over entire range.
We recall that the {\it quantum product}
of two vector fields $v_1$ and $v_2$ is defined by
\[ v_1 \qp v_2 := \gwii{v_1 \, v_2 \, \gua} \, \ga.\]
Define
\begin{equation}
\Psi := \gwiione{E^{2}} + \frac{1}{24} \sum_{\alpha}
     \left<\left< E E \gamma_{\alpha} \gamma^{\alpha} \right>\right>_{0}
    - \frac{1}{2}\sum_{\alpha}  \left(b_{\alpha} (1-b_{\alpha})
         - \frac{b_{1}+1}{6}\right)
    \left<\left< \gamma_{\alpha} \gamma^{\alpha} \right>\right>_{0}
\end{equation}
where $E^2 = E \qp E$ is the quantum square of the Euler vector field.
It was proved in \cite{L1} that, for any smooth projective variety $M$,
the {\it genus-1 Virasoro
conjecture} can be reduced to a single equation on $H^{*}(M; {\Bbb C})$:
\begin{equation} \label{eqn:VirSmall}
    \Psi = 0.
\end{equation}
Moreover, since $E \Psi = \Psi$ (cf. \cite[Lemma 6.3]{L1}),
the genus-1 Virasoro conjecture holds if and only if
\begin{equation}
E \Psi = 0.
\end{equation}
Therefore, to prove the genus-1 Virasoro conjecture, it suffices to show that
$v \Psi =0$ for all vector field $v$ on $H^{*}(M; {\Bbb C})$. It follows from the
string equation that $\gamma_1 \Psi = 0$ where $\gamma_1 = E^0$ is the identity of
the ordinary cohomology ring. In this paper we will give
another vector field which always annihilates $\Psi$.

Define the vector field
\begin{equation}
\Delta := \gua \qp \ga.
\end{equation}
If, in the definition of $\Delta$, we replace the quantum product "$\qp$"
by the ordinary cup product,
we get a vector field proportional to the volume element.
Therefore we call $\Delta$ the {\it quantum volume element}.
The  main result of this paper is the following
\begin{thm} \label{thm:virD}
For all smooth projective varieties,
\[ \Delta \Psi = 0. \]
\end{thm}
This result provides a new evidence for the genus-1 Virasoro conjecture.

\section{Properties of Euler vector fields}
\label{sec:Euler}

We first recall some basic properties of the Euler vector field $E$. We start with the
{\it quasi-homogeneity equation}
\[ \left<\left< E \right>\right>_{g} =
        (3-d)(1-g) F_{g}
        + \frac{1}{2} \delta_{g, 0}\sum_{\alpha, \beta} {\cal C}_{\alpha \beta}
        t^{\alpha}_{0}t^{\beta}_{0}
        - \frac{1}{24} \delta_{g, 1} \int_{M} c_{1}(M) \cup c_{d-1}(M). \]
This equation is a consequence of the divisor equation.
Define the grading operator $G$ by
\[ G(v) := \sum_{\alpha} b_{\alpha} f_{\alpha} \ga \]
for any vector field $v =\sum_{\alpha}  f_{\alpha} \ga$.
Derivatives of quasi-homogeneity equation has the form
\begin{eqnarray}
\gwiig{E \, v_{1} \, \cdots \, v_{k}}
& = & \sum_{i=1}^{k}
    \gwiig{v_{1} \, \cdots \, G(v_{i})
    \, \cdots \, v_{k}}  \nonumber \\
&& - (2g+k-2)(b_{1}+1)
        \gwiig{v_{1} \, \cdots \, v_{k}} \nonumber \\
&& + \delta_{g, 0} \nabla^{k}_{v_{1}, \cdots, v_{k}} \left(
    \frac{1}{2} {\cal C}_{\alpha \beta} t_{0}^{\alpha} t_{0}^{\beta}
    \right) \label{eqn:dhomog}
\end{eqnarray}
where ${\cal C}_{\alpha \beta}$ is defined by $c_1(M) \cup \ga = {\cal C}_{\alpha}^{ \beta} \gb$,
and $\nabla$ is the trivial connection on $H^*(M; \mathbb{C})$ defined by
$\nabla \ga = 0$ for all $\alpha$. In particular,
\begin{equation} \label{eqn:Eg04pt}
 \gwii{E \, v_1 \, v_2 \, \gua} \, \ga
    = G(v_1) \qp v_2 + v_1 \qp G(v_2) - G(v_1 \qp v_2) - b_1 v_1 \qp v_2.
\end{equation}
Combining with \cite[Lemma 4.2]{L1}, we can obtain
\begin{eqnarray}
\nabla_{E^k} \Delta &=& \gwii{E^k \, \gua \, \ga \, \gub} \, \gb \nonumber \\
&=& (k - b_1) E^{k-1} \qp \Delta - G(E^{k-1} \qp \Delta)
    - \sum_{i=1}^{k-1} \, \, \Delta \qp E^{i-1} \qp G(E^{k-i}) \nonumber \\
&&     - \sum_{i=1}^{k-1} G(\Delta \qp E^{i-1}) \qp E^{k-i}
    \label{eqn:dEkD}
\end{eqnarray}
for $k \geq 1$.
Covariant derivative of $E$ is given by
\begin{equation} \label{eqn:dEuler}
 \nabla_{v} E = - G(v) + (b_{1}+1)v.
\end{equation}
Using the fact that
\[ \nabla_w (v_1 \qp v_2) = (\nabla_w v_1) \qp v_2 + v_1 \qp (\nabla_w v_2)
    + \gwii{w \, v_1 \, v_2 \, \gua } \, \ga, \]
we can also show that
\begin{equation} \label{eqn:dDE2}
\nabla_{\Delta} E^2 = \Delta \qp G(E) - G(\Delta) \qp E - G(\Delta \qp E)
    + (b_1+2) \Delta \qp E.
\end{equation}
Combining equations \eqref{eqn:dEkD} and \eqref{eqn:dDE2}, we have
\[ [E^2, \Delta] = - 2 b_1 E \qp \Delta - 2 G(E) \qp \Delta.\]

\section{Proof of the main theorem}

For any vector
fields $v_{1}, \ldots v_{4}$ on the small phase space, we define
\begin{eqnarray*}
G_{0}(v_{1}, v_{2}, v_{3}, v_{4}) & = &
    \sum_{g \in S_{4}} \sum_{\alpha, \beta} \left\{
        \frac{1}{6} \left<\left< v_{g(1)} v_{g(2)}v_{g(3)}
                    \gamma^{\alpha} \right>\right>_{0}
        \left<\left< \gamma_{\alpha} v_{g(4)} \gamma_{\beta}
                    \gamma^{\beta} \right>\right>_{0}
        \right. \\
    && \hspace{40pt}
    + \frac{1}{24} \left<\left< v_{g(1)} v_{g(2)}v_{g(3)} v_{g(4)}
                    \gamma^{\alpha} \right>\right>_{0}
    \left<\left< \gamma_{\alpha}  \gamma_{\beta}
                    \gamma^{\beta} \right>\right>_{0}
            \\
    && \hspace{40pt} \left.
    - \frac{1}{4} \left<\left< v_{g(1)} v_{g(2)}
            \gamma^{\alpha} \gamma^{\beta} \right>\right>_{0}
        \left<\left< \gamma_{\alpha}  \gamma_{\beta}
                v_{g(3)} v_{g(4)} \right>\right>_{0}
            \right\},
\end{eqnarray*}
and
\begin{eqnarray*}
G_{1}(v_{1}, v_{2}, v_{3}, v_{4}) & = &
    \sum_{g \in S_{4}}
        3 \left<\left< \{v_{g(1)} \qp v_{g(2)} \}
            \{v_{g(3)} \qp v_{g(4)} \}
                     \right>\right>_{1}
        \\
    &&
    - \sum_{g \in S_{4}}
    4 \left<\left< \{v_{g(1)} \qp v_{g(2)} \qp
            v_{g(3)} \} v_{g(4)}
                     \right>\right>_{1}
    \\
    &&
     - \sum_{g \in S_{4}} \sum_{\alpha}
        \left<\left< \{ v_{g(1)} \qp v_{g(2)} \}
                v_{g(3)} v_{g(4)}
                    \gamma^{\alpha} \right>\right>_{0}
    \left<\left< \gamma_{\alpha} \right>\right>_{1}
        \\
    &&
     + \sum_{g \in S_{4}} \sum_{\alpha}
         2 \left<\left< v_{g(1)}  v_{g(2)} v_{g(3)}
                    \gamma^{\alpha} \right>\right>_{0}
    \left<\left< \{\gamma_{\alpha} \qp v_{g(4)} \}
            \right>\right>_{1}.
\end{eqnarray*}
Note that $G_{0}$ is completely determined  by genus-0 data, while each term in
$G_{1}$ contains genus-1 information. These two tensors are connected
by Getzler's equation (cf. \cite{Ge}):
\begin{equation} \label{eqn:Ge}
    G_{0} + G_{1} = 0.
\end{equation}
Theorem~\ref{thm:virD} is obtained by applying this equation to
$v_1 = v_2 = E$, $v_3 = \gua$, $v_4 = \ga$, and summing over $\alpha$.

We first consider the genus-$1$ part of Equation \eqref{eqn:Ge}.
\begin{lem} \label{lem:G1EED}
\[ \sum_{\alpha} G_1(E, E, \gua, \ga) = 24 \Delta \gwiione{E^2}. \]
\end{lem}
{\bf Proof}: We will use the convention that repeated indices should be summed
over their entire range. Therefore we will omit $\sum_{\alpha}$ in the left hand of
this formula. To compute $G_1(E, E, \gua, \ga)$, we notice that
\begin{eqnarray*}
 \gwiione{\{ E \qp \ga \} \, \{ \gua \qp E \}}
&=& \gwii{E \ga \gub} \gwii{E \gua \gum} \gwiione{\gb \gm} \\
&=& \gwii{E E \ga} \gwii{\gua \gub \gum} \gwiione{\gb \gm} \\
&=& \gwiione{\{E^2 \qp \gum \} \, \gm}
= \gwiione{\{E^2 \qp \gua \} \, \ga} .
\end{eqnarray*}
In the second equality, we have used the associativity of the quantum product.
This observation enables us to simplify the formula for
$G_1(E, E, \gua, \ga)$ and obtain
\begin{eqnarray}
&& G_1(E, E, \gua, \ga)  \nonumber \\
&=& 24 \gwiione{E^2 \, \Delta} - 48 \gwiione{\{ E \qp \Delta\} \, E}
    - 4 \gwii{E^2 \, \gua \, \ga \, \gub} \gwiione{\gb} \nonumber \\
&&    - 16 \gwii{\{E \qp \gua \} \, \ga \, E \, \gub} \gwiione{\gb}
    - 4 \gwii{\Delta \, E \, E  \, \gub} \gwiione{\gb} \nonumber \\
&&    + 24 \gwii{E \, E \, \gua \, \gub} \gwiione{\{ \ga \qp \gb \}}
    + 24 \gwii{E \, \ga \, \gua \, \gub} \gwiione{\{ \gb \qp E \}}.
        \label{eqn:G1EE}
\end{eqnarray}
We now use formulas in Section~\ref{sec:Euler} to compute each term on the right hand side of this
equation. Using equation~\eqref{eqn:dDE2}, we have
\begin{eqnarray*}
\gwiione{E^2 \, \Delta}
&=& \Delta \gwiione{E^2} - \gwiione{\left\{ \nabla_{\Delta} E^2 \right\}} \\
&=& \Delta \gwiione{E^2} - \gwiione{ \left\{ \Delta \qp G(E) - G(\Delta) \qp E - G(\Delta \qp E)
    + (b_1+2) \Delta \qp E \right\} }.
\end{eqnarray*}
Since $\gwiione{E}$ is a constant due to the quasi-homogeneity equation,
 by equation \eqref{eqn:dEuler},
we have
\begin{eqnarray*}
\gwiione{\{ E \qp \Delta\} \, E}
&=& \{ E \qp \Delta\} \gwiione{ E}  - \gwiione{\left\{ \nabla_{E \qp \Delta} E \right\}} \\
&=& \gwiione{ \left\{ G( E \qp \Delta) - (b_{1}+1)  E \qp \Delta \right\} }.
\end{eqnarray*}
By equation \eqref{eqn:dEkD}, we have
\begin{eqnarray*}
\gwii{E^2 \, \gua \, \ga \, \gub} \gwiione{\gb}
&=& \gwiione{ \left\{ (2-b_1) E \qp \Delta - G(E \qp \Delta) - G(E) \qp \Delta - E \qp G(\Delta)
            \right\} }.
\end{eqnarray*}
By equation \eqref{eqn:Eg04pt}, we have
\begin{eqnarray*}
&& \gwii{\{E \qp \gua \} \, \ga \, E \, \gub} \gwiione{\gb} \\
&=& \gwiione{ \left\{ G(E \qp \gua ) \qp \ga + E \qp \gua  \qp G(\ga)
        - G(E \qp \Delta) - b_1 E \qp \Delta
        \right\} }.
\end{eqnarray*}

As a convention, we arrange the basis
$\{\gamma_{1}, \ldots, \gamma_{N}\}$ of $H^{*}(M, \mathbb{C})$
in such a way  that the degree $p_{\alpha} + q_{\alpha}$
of $\gamma_{\alpha} \in H^{p_{\alpha}, q_{\alpha}}$ is non-decreasing with respect to $\alpha$
and if two cohomology classes have the same dimension, we also require that
the holomorphic dimension $p_{\alpha}$ is non-decreasing. Under this convention, we have
\[ G(\gua) = (1- b_{\alpha}) \gua \]
for all $\alpha$, and
\[ G(\gua) \qp \ga = \Delta - \gua \qp G(\ga).\]
On the other hand,
\[ G(\gua) \qp \ga = \eta^{\alpha \beta} G(\gb) \qp \ga
    = G(\gb) \qp \gub = \gua \qp G(\ga).\]
So we must have
\begin{equation} \label{eqn:Gaua}
 G(\gua) \qp \ga = \gua \qp G(\ga) = \frac{1}{2} \, \Delta.
 \end{equation}
Hence
\begin{eqnarray}
 G(E \qp \gua ) \qp \ga
&=& \gwii{E \, \gua \, \gub} G(\gb) \qp \ga
= G(\gb) \qp (E \qp \gub) \nonumber \\
&=& \frac{1}{2} \, E \qp \Delta. \label{eqn:EGaua}
\end{eqnarray}
Therefore we obtain
\begin{eqnarray*}
 \gwii{\{E \qp \gua \} \, \ga \, E \, \gub} \gwiione{\gb}
&=& \gwiione{ \left\{ (1-b_1) E \qp \Delta
        - G(E \qp \Delta)
        \right\} }.
\end{eqnarray*}
Similarly,
\begin{eqnarray*}
\gwii{\Delta \, E \, E  \, \gub} \gwiione{\gb}
&=& \gwiione{ \left\{
    G(\Delta) \qp E + \Delta \qp G(E) - G(\Delta \qp E) - b_1 \Delta \qp E
    \right\} },
\end{eqnarray*}
and
\begin{eqnarray*}
&& \gwii{E \, E \, \gua \, \gub} \gwiione{\{ \ga \qp \gb \}} \\
&=& \gwiione{ \left\{ \ga \qp \left(
    G(E) \qp \gua + E \qp G(\gua) - G(E \qp \gua) - b_1 E \qp \gua
    \right) \right\} } \\
&=& \gwiione{ \left\{
    G(E) \qp \Delta - b_1 E \qp \Delta
     \right\} }.
\end{eqnarray*}
To compute the last term in equation \eqref{eqn:G1EE}, we first compute
\begin{eqnarray}
  \gwii{E \, \ga \, \gua \, \gub}  \gb
&=&  G(\ga) \qp \gua + \ga \qp G(\gua) - G(\Delta) - b_1 \Delta
   \nonumber  \\
&=& (1-b_1)  \Delta -  G(\Delta) \label{eqn:Eaua}
\end{eqnarray}
by equation \eqref{eqn:Gaua}.
So the last term in equation \eqref{eqn:G1EE} is
\begin{eqnarray*}
 \gwii{E \, \ga \, \gua \, \gub} \gwiione{\{ \gb \qp E \}}
&=& \gwiione{ \left\{
    (1-b_1) E \qp \Delta - E \qp G(\Delta)
    \right\} }.
\end{eqnarray*}
After plugging the above formulas into equation \eqref{eqn:G1EE}, all terms on the
right hand side cancel except the term $24 \Delta \gwiione{E^2}$.
The lemma is thus proved.
$\Box$

Now we consider the genus-$0$ part of Equation \eqref{eqn:Ge}.
Let
\begin{equation} \label{eqn:PhiDef}
\Phi := - \frac{1}{24} \sum_{\alpha}
     \left<\left< E E \gamma_{\alpha} \gamma^{\alpha} \right>\right>_{0}
    + \frac{1}{2}\sum_{\alpha}  \left(b_{\alpha} (1-b_{\alpha})
         - \frac{b_{1}+1}{6}\right)
    \left<\left< \gamma_{\alpha} \gamma^{\alpha} \right>\right>_{0}.
\end{equation}
Then
\begin{equation}
 \Psi = \gwiione{E^2} - \Phi
\end{equation}
and the genus-1 Virasoro conjecture can be reduced to
\[ \gwiione{E^2} = \Phi. \]
\begin{lem} \label{lem:G0EED}
\[ \sum_{\alpha} G_0(E, E, \gua, \ga) = -24 \Delta \Phi. \]
\end{lem}
{\bf Proof}:
Again we will assume that repeated indices will be summed over their entire range.
First, by definition of $G_0$, we have
\begin{eqnarray}
G_0(E, E, \gua, \ga)
&=& 2\gwii{E \, E \, \gua \, \gub}\gwii{\gb \, \ga \, \gum \, \gm}
    + 2 \gwii{E \, \gua \, \ga \, \gub}\gwii{\gb \, E \, \gum \, \gm}
    \nonumber \\
&&   + \gwii{E \, E \, \gua \, \ga \, \Delta}
    - 2 \gwii{E \, E \, \gub \, \gum} \gwii{\gb \, \gm \, \gua \, \ga}
    \nonumber \\
&&    - 4 \gwii{E \, \gua \, \gub \, \gum} \gwii{\gb \, \gm \, E \, \ga}.
  \label{eqn:G0EE}
\end{eqnarray}
Note that the first and the fourth terms on the right hand side are canceled
with each other.
Applying equation \eqref{eqn:Eaua} to the second term, we have
\begin{eqnarray}
 \gwii{E \, \gua \, \ga \, \gub}\gwii{\gb \, E \, \gum \, \gm}
&=& \gwii{\left\{ (1- b_1) \Delta - G(\Delta) \right\}
            \, E \, \gum \, \gm}. \label{eqn:EDED}
\end{eqnarray}
Using equation \eqref{eqn:Eaua} again, we obtain
\begin{eqnarray}
\gwii{ \Delta \, E \, \gum \, \gm}
&=& \gwii{E \, \gum \, \gm \, \gub} \gwii{\gb \, \gua \, \ga} \nonumber \\
&=& \gwii{ \left\{(1-b_1) \Delta - G(\Delta) \right\} \gua \, \ga}. \nonumber
\end{eqnarray}
Moreover,
\begin{eqnarray*}
\gwii{G(\Delta) \, \gua \, \ga}
&=& b_{\mu} \gwii{\gub \, \gb \, \gum} \gwii{\gm \, \gua \, \ga} \\
&=& \gwii{\gub \, \gb \, \left\{\gum - G(\gum) \right\}} \gwii{\gm \, \gua \, \ga} \\
&=& \gwii{\Delta \, \gua \, \ga} - \gwii{\gub \, \gb \, G(\Delta)}.
\end{eqnarray*}
Moving the second term on the right hand to the left hand, we obtain
\begin{equation} \label{eqn:GDD}
\gwii{G(\Delta) \, \gua \, \ga} = \frac{1}{2} \, \gwii{\Delta \, \gua \, \ga}.
\end{equation}
Hence, we have
\begin{eqnarray}
\gwii{ \Delta \, E \, \gum \, \gm}
&=& \left(\frac{1}{2} - b_1 \right) \gwii{\Delta \, \gua \, \ga}. \label{eqn:EDD}
\end{eqnarray}
By equation \eqref{eqn:EDED}, we have
\begin{eqnarray}
&& \gwii{E \, \gua \, \ga \, \gub}\gwii{\gb \, E \, \gum \, \gm} \nonumber \\
&=& (1- b_1) \left(\frac{1}{2} - b_1 \right) \gwii{  \Delta \, \gua \, \ga}
    - \gwii{ G(\Delta) \, E \, \gum \, \gm}. \label{eqn:EDED2}
\end{eqnarray}

To compute the last term on the right hand side
of equation \eqref{eqn:G0EE}, we set
\[ f := \gwii{E \, \gua \, \gub \, \gum} \gwii{\gb \, \gm \, E \, \ga}. \]
Applying equation \eqref{eqn:Eg04pt}, we have
\begin{eqnarray}
f
&=&  \gwii{\left\{ G(\gua) \qp \gum + \gua \qp G(\gum) - G(\gua \qp \gum)
        - b_1 \gua \qp \gum \right\} \, \gm \, E \, \ga} \nonumber \\
&=& (2-b_{\alpha}- b_{\mu}-b_1) \gwii{\left\{ \gua \qp \gum \right\} \, \gm \, E \, \ga}
    -  \gwii{ G(\gua \qp \gum) \, \gm \, E \, \ga} \nonumber \\
&=& (2-2b_{\alpha}-b_1) \gwii{\gua \, \gum \, \gub} \gwii{\gb \, \gm \, E \, \ga}
    - \gwii{\gua \, \gum \, \gub} \gwii{ G(\gb) \, \gm \, E \, \ga}. \nonumber
\end{eqnarray}
Switching $\alpha$ and $\beta$ in the last term, we have
\begin{equation}
f = (2-b_1)\gwii{E \, \ga \, \gb \, \gum}\gwii{\gm \, \gua \, \gub}
    -3\gwii{E \, G(\ga) \, \gb \, \gum} \gwii{\gm \, \gua \, \gub}.
    \label{eqn:f3}
\end{equation}
Applying equation \eqref{eqn:Eg04pt} again, we have
\begin{eqnarray}
&& \gwii{E \, \ga \, \gb \, \gum}\gwii{\gm \, \gua \, \gub} \nonumber \\
&=& \gwii{ \left\{ G(\ga) \qp \gb + \ga \qp G(\gb) - G(\ga \qp \gb)
        - b_1 \ga \qp \gb \right\} \, \gua \, \gub}.
            \nonumber
\end{eqnarray}
By the associativity of the quantum product and equation \eqref{eqn:Gaua},
\begin{eqnarray}
\gwii{ \left\{ G(\ga) \qp \gb \right \} \, \gua \, \gub}
&=& \gwii{ \left\{ G(\ga) \qp \gua \right \} \, \gb \, \gub}
    \,
= \, \frac{1}{2} \, \gwii{ \Delta \, \gua \, \ga}
    \label{eqn:GD}
\end{eqnarray}
and
\begin{eqnarray}
\gwii{  G(\ga \qp \gb)  \, \gua \, \gub}
&=& \gwii{\ga \, \gb \, \gum} \gwii{  G(\gm)  \, \gua \, \gub}
    \nonumber \\
&=& \gwii{ \ga \, \left\{ G(\gm) \qp \gua \right\} \, \gum}
= \frac{1}{2} \, \gwii{ \Delta \, \gua \, \ga}.
\end{eqnarray}
So we have
\begin{eqnarray}
 \gwii{E \, \ga \, \gb \, \gum}\gwii{\gm \, \gua \, \gub}
&=& \left(\frac{1}{2} - b_1 \right) \gwii{\Delta \, \gua \, \ga}.
\label{eqn:EDDD}
\end{eqnarray}
Moreover
\begin{eqnarray}
&& \gwii{E \, G(\ga) \, \gb \, \gum} \gwii{\gm \, \gua \, \gub} \nonumber \\
&=& \gwii{ \left\{ G(G(\ga)) \qp \gb + G(\ga) \qp G(\gb)
        - G(G(\ga) \qp \gb) - b_{1} G(\ga) \qp \gb \right\} \, \gua \, \gub}.
        \nonumber
\end{eqnarray}
Since
\begin{eqnarray}
&& \gwii{\left\{ G(\ga) \qp G(\gb) \right\} \, \gua \, \gub} \nonumber \\
&=& \gwii{\left\{ G(\ga) \qp  \gua \right\} \, G(\gb) \, \gub}
= \frac{1}{2} \, \gwii{\left\{ \ga \qp  \gua \right\} \, G(\gb) \, \gub} \nonumber \\
&=& \frac{1}{2} \, \gwii{\ga \,  \gua  \, \left\{ G(\gb) \qp \gub \right\} }
= \frac{1}{4} \, \gwii{ \Delta \, \gua \, \ga}
\end{eqnarray}
and
\begin{eqnarray}
&& \gwii{G( G(\ga) \qp \gb ) \, \gua \, \gub} \nonumber \\
&=& \gwii{G(\ga) \, \gb \, \gum} \gwii{G(\gm) \, \gua  \, \gub}
=  \gwii{\left\{ G(\ga) \qp  \gum \right\} \, G(\gm) \, \gua} \nonumber \\
&=& \gwii{\left\{ G(\ga) \qp  \gua \right\} \, G(\gm) \, \gum}
= \frac{1}{2} \gwii{\left\{ \ga \qp  \gua \right\} \, G(\gm) \, \gum}
    \nonumber \\
&=& \frac{1}{2} \gwii{ \ga \,  \gua \, \left\{ G(\gm) \qp \gum \right\} }
= \frac{1}{4} \, \gwii{ \Delta \, \gua \, \ga},
\end{eqnarray}
together with equation \eqref{eqn:GD}, we have
\begin{eqnarray}
 \gwii{E \, G(\ga) \, \gb \, \gum} \gwii{\gm \, \gua \, \gub}
&=& \left(b_{\alpha}^2 - \frac{1}{2} b_{1} \right) \gwii{\Delta \, \gua \, \ga}.
\label{eqn:EGD}
\end{eqnarray}
Combining results of equations \eqref{eqn:f3}, \eqref{eqn:EDDD}, and \eqref{eqn:EGD}, we obtain
that the last term on the right hand side
of equation \eqref{eqn:G0EE} is
\begin{equation}
f = (-3 b_{\alpha}^2 + b_{1}^2 - b_{1} + 1) \gwii{\Delta \, \gua \, \ga}.
\end{equation}
Together with equation \eqref{eqn:EDED2}, we can simplify equation \eqref{eqn:G0EE} as
\begin{eqnarray}
G_0(E, E, \gua, \ga)
&=& \gwii{E \, E \, \gua \, \ga \, \Delta}
    - 2 \gwii{ G(\Delta) \, E \, \gum \, \gm} \nonumber \\
&&  + (12 b_{\alpha}^2 - 2 b_{1}^2 + b_{1} -3) \gwii{\Delta \, \gua \, \ga}.
  \label{eqn:G0EE2}
\end{eqnarray}

On the other hand, by the definition of $\Phi$ in equation \eqref{eqn:PhiDef}, we have
\begin{eqnarray}
24 \Delta \Phi
&=& - \gwii{ \Delta \, E \, E \, \ga \, \gua }
    - 2 \gwii{ \left\{ \nabla_{\Delta} E \right\} \, E \, \ga \, \gua }
    \nonumber \\
&&    + 12   \left(b_{\alpha} (1-b_{\alpha})
         - \frac{b_{1}+1}{6}\right)
    \gwii{ \Delta \, \ga \, \gua}. \nonumber
\end{eqnarray}
By equations \eqref{eqn:dEuler} and \eqref{eqn:EDD}, we have
\begin{eqnarray*}
\gwii{ \left\{ \nabla_{\Delta} E \right\} \, E \, \ga \, \gua }
&=& \gwii{ \left\{ - G(\Delta) + (b_{1}+1) \Delta \right\} \, E \, \ga \, \gua } \\
&=& - \gwii{ G(\Delta) \, E \, \ga \, \gua }
    + (b_{1}+1) \left(\frac{1}{2} - b_1 \right) \gwii{\Delta \, \gua \, \ga}.
\end{eqnarray*}
Moreover,
\[ b_{\alpha}\gwii{\Delta \, \gua \, \ga} = \gwii{\Delta \, \gua \, G(\ga)}
    = \frac{1}{2} \gwii{\Delta \, \gua \, \ga}.\]
So we have
\begin{eqnarray}
24 \Delta \Phi
&=& - \gwii{ \Delta \, E \, E \, \ga \, \gua }
    + 2  \gwii{ G(\Delta) \, E \, \ga \, \gua } \nonumber \\
&&  + \left(- 12 b_{\alpha}^2 + 2 b_1^2 - b_1 + 3 \right)
    \gwii{ \Delta \, \ga \, \gua}. \nonumber
\end{eqnarray}
Comparing with equation \eqref{eqn:G0EE2}, we obtain
\[ G_0(E, E, \gua, \ga) = - 24 \Delta \Phi. \]
The lemma is thus proved.
$\Box$

\vspace{10pt}
\noindent
{\bf Proof of Theorem~\ref{thm:virD}}:
Since $\Psi= \gwiione{E^2}- \Phi$, this theorem follows from
Lemmas \ref{lem:G1EED}, \ref{lem:G0EED} and Equation \eqref{eqn:Ge}.
$\Box$



\vspace{20pt}
 \noindent
Xiaobo Liu \\

\noindent
Beijing International Center for Mathematical Research, \\
Beijing University, Beijing, China. \\
\& \\
Department of Mathematics, \\
University of Notre Dame, \\
Notre Dame, IN 46556, USA \\

\noindent
Email: {\it xliu3@nd.edu}

\end{document}